\documentclass{amsart}

\newcommand{\Fin}{F_{\infty}}
\newcommand{\sub}{\subseteq}
\newcommand{\cont}{\supseteq}

\newcommand{\normal}{\triangleleft}
\newcommand{\W}{{\mathfrak W}}
\newcommand{\V}{{\mathfrak V}}
\newcommand{\Ni}{{\mathfrak N}}
\renewcommand{\O}{{\mathfrak O}}
\newcommand{\U}{{\mathfrak U}}
\newcommand{\A}{{\mathfrak A}}
\newcommand{\N}{{\mathbb N}}
\newcommand{\Z}{{\mathbb Z}}
\newcommand{\Zs}{{\mathbb Z}}
\newcommand{\Q}{{\mathbb Q}}

\newcommand{\sn}{{\normal  \normal \, }}

\newcommand{\Inn}[1]{\mathrm{Inn}\left( #1 \right)}
\newcommand{\Aut}[1]{\mathrm{Aut}\left( #1 \right)}
\newcommand{\card}[1]{\mathrm{card}\left( #1 \right)}

\newcommand{\var}[1]{\mathrm{var}\left( #1 \right)}

\newcommand{\Diag}[1]{\mathrm{Diag}\left( #1 \right)}
\newcommand{\Wr}{\,\mathrm{Wr}\,}

\newcommand{\proofend}{{\hfill $\lefteqn{\sqcap}\sqcup$}
\vskip4mm}

\theoremstyle{definition}
  
\theoremstyle{plain}
  
\newtheorem{Theorem}{\sc Theorem}  
\newtheorem{Proposition}{\sc Proposition}

\begin{document}

\renewcommand{\subjclassname}{2000 Mathematics Subject Classification}

\title{Subnormal embedding theorems for groups\footnote{\vskip1mm \centerline{Journal of the London Mathematical Society, 62.2 (2000), 398--406.}}}

\author{Vahagn H. Mikaelian}

\date{30 May 1998}

\dedicatory{\small To the memory of my brother Sassoun Mikaelian}  

\begin{abstract}
In this paper we establish some subnormal embeddings of groups into groups with additional properties; in particular embeddings of countable groups into $2$-generated groups with some extra properties. The results obtained are generalizations of theorems of P.~Hall, R.~Dark, B.~Neumann, Hanna Neumann, G.~Higman on embeddings of that type. Considering subnormal embeddings of finite groups into finite groups with additional properties we bring some illustration to a result of H.~Heineken and J.~Lennox.
\end{abstract}

\maketitle


\section{Introduction and the main results}
\label{Introductionandmainstatements}

The aim of this paper is to introduce a method of  subnormal ``verbal"
embeddings of groups into larger groups with additional properties. 
First we prove 
the following main result:

\begin{Theorem}
\label{Interesting}

Let $H$ be an arbitrary infinite group and $V \sub \Fin$ be an 
arbitrary non-trivial word set. Then
 
 {\bf A.} there exists a group $G$ of the same cardinality as $H$ 
 such that $G$ has a 
 subnormal subgroup $\tilde H \cong H$ and $\tilde H \sub V(G)$;
 \vskip1mm
 {\bf B.} if the group $H$ is countable the (countable) group $G$      
 can be chosen two-generated; 
\vskip1mm
 {\bf C.} if the group $H$ is solvable the group $G$ of the statement 
 {\bf A} or  the two-generated group  of the statement {\bf B} can be 
 chosen as  well solvable.
   
\end{Theorem} 
  
This theorem is also true for finite groups, as  a special case of countable 
groups. However not each {\sl finite} group $H$  can be subnormally 
embedded into some {\sl finite} group $G$ such that the isomorphic image 
of $H$ lies in $V(G)$. We will return to 
the case of finite groups later. 

To avoid repetition we will sometimes talk of the group
$H$ and of its  image $\tilde H$ in $G$ as of the same group.

In fact the statement {\bf A} establishes a subnormal embedding of 
defect two. 
Thus it is natural to ask whether this defect  can be reduced to one, i. e. 
whether the ``verbal'' embedding of $H$ in $G$ can be {\sl normal}. 
The answer is contained in Theorem~\ref{CombinationWithAnOther} in Section 3.

The statement {\bf B} is a generalization  of a theorem of 
G.~Higman, B.~Neumann and 
Hanna Neumann on embeddability of an arbitrary countable group in 
a $2$-generated group
(1947,~\cite{HigmanNeumannNeumann}). The theorem, 
called in the book of 
Robinson~\cite{Robinson}
``probably the most famous of all embedding theorems'', was a stimulus for
further research on embeddings into $2$-generated groups and related
problems. See the papers of  Ol'shanskii~\cite{Ol'shanskiiVestnikMoskov} 
\cite{Ol'shanskiiUkrainMatZh}, Galvin~\cite{Galvin}, 
B. Neumenn~\cite{NeumannNieuwArch},   Wilson and 
Zalesskii~\cite{WilsonZalesskii},
Levin and Rosenberger~\cite{LevinRosenberger},
and literature cited there.

The idea that this 
embedding can be  subnormal belongs to P.~Hall and R.~Dark~\cite{Dark}. In
fact the Hall--Dark construction can be slightly modified 
to get a  subnormal embedding
into the {\sl derived group} of some two-generated group. 

The statement {\bf C}  generalizes the result of B. Neumann and Hanna 
Neumann~\cite{NeumannHNeumann} on embeddability of {\sl solvable} 
countable groups in the second derived groups of  {\sl solvable} 
two-generated groups. Considering in some sense more ``economical'' embeddings 
B.~Neumann and Hanna Neumann proved that the solvable group $H$ 
of length $n$ can be embedded into a
two-generated solvable group $G$  of length $n+2$ but not in 
general of length $n+1$, this is  Corollary 5.2 of Theorem 
5.1~\cite{NeumannHNeumann}    saying that each group
of variety $\U$ can be embedded into some two-generated  group 
of variety $\U \A^2$, where 
$\A$ is the variety of all abelian groups. This result cannot 
be strenghtened by the requirement to embed  into a verbal 
subgroup (see below).   Nevertheless a simplified version 
of our construction can be used in order to make both embeddings 
of the mentioned Theorem 5.1 and Corollary 5.2~\cite{NeumannHNeumann} 
{\sl subnormal} (respectively Theorems~\ref{VtimesAA}
and~\ref{SolvableIntoSolvable} in Section~\ref{OnSolvables}). 

The importance of the case of solvable groups follows, in particular, from
the fact that we could {\sl not} obtain the  analogs of the mentioned 
result for {\sl abelian} 
or {\sl nilpotent} groups because each subgroup of a two-generated 
abelian group  should be as well  two-generated, and because each 
subgroup of a finitely generated nilpotent
group  should be as well finitely generated  \cite{Robinson}. The countable 
abelian group, say, $\Q$ of rational numbers is infinitely 
generated and cannot be contained in any
finitely generated nilpotent group. 
\vskip2mm

The statement {\bf A} of   Theorem~\ref{Interesting} is 
in contrast with 
the {\sl negative} results regarding  subnormal embeddings 
of {\sl finite} groups
into {\sl finite} groups. P.~Hall and R.~Dark have shown that for the 
symmetric group
$S_3$ there is no finite $G$ such that 
$S_3 \sn G$ and $S_3 \sub G'$ (see 
in~\cite{Dark}). This result 
is strongly generalized by 
Heineken and Lennox~\cite{HeinekenLennox}
who proved that for {\sl no} finite solvable 
complete group $H$ there exists a finite $G$
such that $H \sn G$ and $H \sub G'$ (see also the extension of
the latter result by Heineken and Soules~\cite{HeinekenSoules}). 

We make this contrast more visible in Section 5
and show that the result of Heineken and Lennox
finds almost the ``smallest'' possible finite groups of that type.
Namely we prove first that for each finite {\sl nilpotent} group $H$
and each non-trivial set of 
words $V$ there exists a finite (and nilpotent)
group $G$ such that $H \sn G$ and $H \sub V(G)$  
(Proposition~\ref{FiniteNilpotentIntoFiniteNilpotent}).

Thus it is very natural to ask 
whether there exist non-nilpotent (non-solvable)
finite groups $H$ for which 
it is possible (impossible) to find a finite 
group $G$ such that $H \sn G$ and $H \sub V(G)$. 
Proposition~\ref{VeryManyGroups} 
says something more general; for {\sl each}
variety $\U$ different from the variety $\O$ of all groups 
there exist finite groups 
{\sl not} from $\U$ that have (or do not have)
embeddings of the mentioned type. And it is possible to 
find such groups in some very small classes of groups 
(Proposition~\ref{ManyGroupsInSmallClasses}). 

This shows
that the question of existence of subnormal ``verbal'' 
embeddings of finite groups
into finite groups is far from being a trivial task. We dare formulate it as a
{\sl problem}  at the end of the current paper.

\vskip5mm
\begin{center}
{\bf Acknowledgments}
\end{center}

It is a great pleasure for the author
to express his gratitude to Professor Hermann 
Heineken for his attention to my research project and for  useful 
considerations while I enjoyed the very warm hospitality of the Mathematics
Institute of the W\"urzburg University, Germany.
The author is thankful to the German Academic 
Exchange Service for financial support (Grant Nr. A/97/13683).


  \section{The main construction}
  \label{Themainconstruction}
  
The following steps prove   Theorem~\ref{Interesting}. 

\vskip3mm

{\bf a.} Let $L$ be a nilpotent 
countable group such that its verbal subgroup 
$V(L)$ contains a subgroup 
isomorphic to $\Z$. To construct $L$ we note
that the set of all nilpotent 
groups generates the variety $\O$ of all groups 
(free groups are residually finite $p$-groups). The variety
$$
                \V = \var {\Fin / V(\Fin) }
$$
corresponding to the word set $V$ is different from $\O$.
 Thus there exists a nilpotent group $L= F_n(\Ni_c)$
of some class $c$    and some  rank $n$ such that $L \notin \V$. An 
element $a$ of infinite order in non-trivial 
and torsion-free subgroup $V(L)$ 
generates an isomorphic copy of  
$\Z =\langle 1_{\Zs} \rangle \cong \langle a \rangle$. 
 $1_{\Zs}$ is not the trivial element of $\Z$ (unlike in the 
 notions that will come later) but the integer $1$.
 
\vskip3mm 
 
{\bf b.} In order to define a 
special linear  order $\prec$ on the countable group 
$L$ we establish the one-to-one
correspondence $\alpha$ between the additive group 
of integers $\Z = \langle 1_{\Zs} \rangle$ and $L$:
$$
           \alpha ( i)= \alpha ( i\cdot 1_{\Zs} ) \,\,\, \hbox{{\rm is equal to}}\,
                  \left\{
                  \begin{array}{l}
                  a^i \,\,\, \hbox{{\rm if and only if}}\,\,\, 
                  i \in \{ \ldots , -2, -1, 0\},\\
                  \hbox{{\rm  an arbitrary $x \in  
                  L \backslash \{a^j\}_{j \le 0} $ 
                  if $i \in \{1,2,\ldots\}$}} 
                  \end{array}
                  \right.
$$
in such a manner that  $\alpha$ is, in particular, a bijection between 
the sets
$$
\Z \backslash \{ \ldots , -2, -1, 0\} \quad \mathrm{and} \quad 
L \backslash \{   \ldots , a^{-2}, a^{-1}, 1_L \}.
$$
\vskip2mm
Now define
$$
\forall \,\, l_1, l_2 \in L, \quad \! l_1\prec l_2 \quad \! 
\hbox{{\rm if and only if}} \quad \! \alpha^{-1}(l_1)\le\alpha^{-1}(l_2).
$$

\vskip3mm

{\bf c.} Consider the 
(cartesian) wreath product $H\Wr L$ and 
assume that the copies of $H$ in 
the base group $H^{(L)}$  are also ordered
according to $\prec$. To be brief 
we write simply $h_i$ instead of $h_{\alpha(i)}$ 
for element contained in $\alpha(i)$'th copy of $H$. Now 
we define 
the sets $D_k$ ($k \in \N\cup \{ 0\}$) of elements of
the base group as follows:
$
\overline h = (\ldots, h_{-1}, h_{o}, h_{1}, \ldots)  
 \in  D_k$ 
if and only if there exist 
integers $n_1, n_2, \ldots,n_k$  such that 
$n_i < n_{i+1}$, $i=1,\ldots, k-1 $ and

$$
           h_i=
                  \left\{
                  \begin{array}{ll}
                  h_{n_1} & {\rm if} \,\,\,   i\le n_1, \\
                  h_{n_2} & {\rm if} \,\,\,   n_1 <i\le n_2, \\
                  \,\,\, \vdots  &  \quad \quad\vdots  \\
                  h_{n_k} & {\rm if} \,\,\,  n_{k-1} <i\le n_k,\\
                  h_{n_k+1} & {\rm if} \,\,\,   n_{k} <i.             
                 \end{array}
                  \right.
$$
Clearly $\card{D_k}=\card{H}$, because for each $\overline h = 
(\ldots, h_{-1}, h_{o}, h_{1}, \ldots) \in D_k$ 
we have only finitely many ``jump" point 
indices $n_1, n_2, \ldots,n_k$; for
each such an index  one can find only countably many ``places'' 
and because
for each index $n_i$  ($i=1, \ldots, k$) and for index $n_k+1$ 
the values of $h_i$
have only $\card{H}$ possibilities.

The group $L$ is countable. Therefore the subgroup 
$$
G_1= \langle L,\bigcup_{i=1}^{\infty} D_k \rangle
$$
of  $H\Wr L$ is of cardinality $\card H$ for it is generated 
by an infinite set 
of cardinality $\card H$.

\vskip3mm

{\bf d.}  The group $H$ is isomorphic to its first copy in the 
base group $H^{(L)}$. This copy lies in $G_1$ because
$$
\overline h_0 = (\ldots, 1_H, h, 1_H, \ldots) \in D_2  
$$
(the element $h$ has the number $0$).
Let us further take for an arbitrary 
$h\in H$ an element $\theta(h) \in H^{(L)}$:
$$
     \theta (h) =(\ldots, h^{-1}, h^{-1}, h^{-1},  
     1_H,  1_H,  1_H, \ldots)\in D_1 \sub G_1
$$
(the first element $1_H$ stands in the place number $0$).
Since 
$$
       a \in V(L) \sub V(G_1)
$$
and since verbal subgroups of groups are characteristic, 
we assert that the
element
$$    a^{\theta(h)} =a\cdot(\ldots,  h, h, 1_H,  
                  1_H, \ldots)^a \cdot (\ldots,  h^{-1}, h^{-1},  
                  1_H,  1_H,  \ldots) 
$$
belongs to $V(G_1)$. The order relation $\prec$ is defined
over $L$ in such a manner that the elements standing in the places
number $\ldots, -2, -1, 0$ go under operation of 
$a$ (right regular representation) one step ``higher'' and
the elements standing in the places number $1, 2, 3, \ldots$ go to
some new places with  numbers 
again from the set $\{1, 2, 3, \ldots \}$. 
Thus
$$
(\ldots,  h, h, 1_H,  1_H, \ldots)^a =(\ldots,  h, h, h,  1_H, \ldots),
$$
where in the second string the first element $1_H$ stands now	 in the 
place number $1$. We get
$$
       a^{\theta(h)} =a\cdot (\ldots, 1_H, h, 1_H, \ldots) 
                  =a\cdot\overline{h}_0 \in V(G_1).
$$
Thus $\overline h_0 \in V(G_1)$ for each $h\in H$. This proves the 
statement {\bf A} of Theorem~\ref{Interesting} 
because the first copy
of $H$ in the base group is clearly subnormal in wreath product 
$H\Wr L$ and in $G=G_1$. And it is clear that if the group $H$ is 
solvable, the group $G\sub H\Wr L$ and its subgroup $G$ are 
solvable, too.

\vskip3mm

{\bf e.} Next we need a group $M$ such that $\card M = \card{G_1}$
and $G_1 \sn M'$. We can simply use the statement {\bf A} we just
proved. But for this case of the {\sl commutator} word $V$ we could 
build  a much smaller construction.  That will
be done in the proof of Theorem~\ref{VtimesAA}.

\vskip3mm

{\bf f.} Finally we use an idea of P. Hall from~\cite{PHall'sIdea}. 
 Let us  $H$ be a countable group subnormally embedded into
the derived subgroup of the 
group $M$ and assume that the elements of the 
countable group $M$ are linearly ordered:
$$
       M=\{ m(0),  m(1),  m(2),   \ldots   \}.
$$
We define:
$$
    m_i =           \left\{
 \begin{array}{cl}
 m(s)  & \hbox{{\rm if $i=2^s$, for some
 $s\in \{ 0, 1, 2,  \},$ }}\\
 1_M   & \hbox{{\rm if $i\not=2^s$, 
 $\hbox{\rm for arbitrary}\, s \in \{ 0, 1, 2,  \}$  }}
 \end{array}
 \right.
$$
and 
$$
  \overline{m}= (\ldots,  m_{-1}, m_{0},m_{1},\ldots )
  \in M^{(\Zs)}\sub M \Wr \Z.
$$
Then for the intersection $T =\langle -1_{\Zs}, 
\overline{m}\, \rangle \cap M^{(\Zs)}$ 
we have according~\cite{PHall'sIdea}:
$$
   T'=  \overline{M'}^{(\Zs)}=K , 
$$
where the latter means {\sl direct} product of copies of $M'$.  

$K$ is normal in $G_2 = \langle -1_{\Zs}, \overline{m} \rangle$ for
it is normal even in $M \Wr \Z$. 
Therefore the subnormal embedding we
are looking for can be given as:
$$
    H \sn G_1 \sn M'  \normal K \normal G_2.
$$
(As we mentioned above, in order to avoid huge formulas, we
talk of a group and its isomorphic image in the larger group 
as of the same group.) 
This proves the statement {\bf B} 
of Theorem~\ref{Interesting}. And the
part of the statement {\bf C} 
concerning embeddability into two-generated
solvable groups follows from the fact that all  groups evolved in
points {\bf e} and {\bf f} above are solvable, too.

Theorem~\ref{Interesting} is proved. \proofend


  \section{The defects of  embeddings,\\
   the case of normal embeddability}
  \label{Thedefectsofbuiltembeddings}

The embedding of $H \sn G_1$ constructed in points {\bf a},...,{\bf d}
of Section~\ref{Themainconstruction} is of defect $2$. The first copy 
of $H$ is subnormal of defect $2$ in $H \Wr L$ and $G_1$ is a 
subgroup of $H \Wr L$ containing this first copy. 

A natural question which arises is: whether
this defect {\sl can} (or {\sl cannot})
be made smaller (namely $1$), i. e. 
\begin{quote}
Do there exist groups $H$
that can (or cannot) be {\sl normally} 
embedded into appropriate groups $G$ such that
the isomorphic image $\tilde H$ of $H$ lies in $V(G)$ for the
given non-trivial word set $V$?
\end{quote}
The general answer to this question (not only for the case of
infinite groups) is given 
in our co-work with Heineken~\cite{HeinekenMikaelianOnnVEmb}, we can assert:
\begin{Theorem}
\label{CombinationWithAnOther} 
Let us $H$ be an arbitrary infinite group and $V \sub \Fin$ be an 
arbitrary non-trivial
word set. Then there exists a group $G$ of the same cardinality 
as $H$ such that $G$ has a 
subnormal subgroup $\tilde H \cong H$ of defect at most $2$ and 
$\tilde H \sub V(G)$. This embedding is of defect $1$ (i. e. normal)
if and only if $V(\Aut H)\cont \Inn H$.
\end{Theorem} 
\begin{proof}
This   follows from Theorem~\ref{Interesting}, from the  
remark made above, from Theorem~1 
of~\cite{HeinekenMikaelianOnnVEmb} and from the fact that if
the infinite group $H$ is normally embeddable into some group $G$ and lies 
in $V(G)$, then  without loss of generality  $G$ is of the same cardinality as
$H$. Let us show it.

The infinite group 	$H$ has a system of generators of the same cardinality as
$H$. So the image $\tilde H \sub V(G)$ is  generated
by the appropriate images, say,  $\tilde h_i$, $i\in I$. For
each element $\tilde h_i \in V(G)$  there is a representation
$$
  \tilde h_i = \left(w_1^{(i)}(g_{11}^{(i)}, 
  \ldots, g_{1q_1}^{(i)})\right)^{\delta_1^{(i)}}\!\!\!\cdots \,\,
  \left(w_u^{(i)}(g_{u1}^{(i)}, \ldots, 
  g_{uq_u}^{(i)})\right)^{\delta_u^{(i)}}\!\!\!\!\!\!,
$$
where $w_1^{(i)}\!\!\!,...\, , w_u^{(i)}$ are some words from $V$, elements 
$\delta_u^{(i)}$ take values $1$ or $-1$, and where $g_{jk}^{(i)}\in G$.
Let us $G_1$ be the subgroup of $G$  generated by this   set $Q$ of elements
$g_{jk}^{(i)}$.
Of course $\tilde H \sub V(G_1)$ and $\tilde H \normal G_1$. Furthermore
$Q$ is of the same cardinality as $H$, and so is $G_1$. \end{proof}

For concrete  examples of groups that do not have normal 
embeddings of mentioned type see~\cite{HeinekenMikaelianOnnVEmb}.

The obtained embedding of countable groups into  two-generator groups
has a larger defect, namely $6$  (because the defect of
$H$ in $G_1$ is  $2$ and because the 
defect of the latter in $M'$ is also
$2$).   

Computation  of an other nature can reduce this defect to $4$. 
But we will not include that here since that would demand an aparatus
strongly different from our main construction.


  \section{The ``economical'' embeddings\\
   of solvable groups into solvable groups}
  \label{OnSolvables}

The initial step here  are the results of B. Neumann and Hanna Neumann, 
mentioned in the introduction,   on embedding of a countable group
of a variety $\U$  into some two-generated  group of the 
variety $\U \A^2$ and
on embedding of a countable solvable group of length  $l$ into a 
two-generated solvable group of length  at most $l+2$.
(Theorem 5.1 and Corollary 5.2 in~\cite{NeumannHNeumann})
In fact these both embeddings are embeddings not only
into the appropriate two-generated group $G$ but into the
{\sl second derived group} of $G$.

Mentioned results are in some sense the best possible ones, since some 
solvable groups of length  $l$ {\sl cannot} be embedded into 
two-generated solvable groups of length  $l+1$. As an example see
Lemma 5.3 in~\cite{NeumannHNeumann}. 

As we mentioned above we cannot add the condition with 
verbal subgroups to
these results. For 
if we take as $V$, say, the identity $\delta_l 
(x_1, \dots , x_{2^l})$ of solvability
of length  $l$: $\delta_0 = x$ and
$$  \delta_{l+1}(x_1, \dots , x_{2^{l+1}}) =   
[\delta_l(x_1, \dots , x_{2^l}), 
  \delta_l(x_{2^l +1}, \dots , x_{2^{l+1}})]
$$
with $l \ge 3$, we get
that for all $n \in \N$ the verbal subgroup $\delta_l 
(G)$ of each solvable group
$G$ of length $n+2$ is of length $n+2-l< n$. 

But we are able to use an easier version of our 
construction to make the 
embeddings of the mentioned results {\sl subnormal}. Namely we build
the construction mentioned in point {\bf e} of 
Section~\ref{Introductionandmainstatements} and 
prove:
\begin{Theorem}
\label{VtimesAA}
For an arbitrary countable group $H$ of the variety $\U$
there exists in the variety $\U \A^2$ a two-generated  group $G$ 
with a subnormal subgroup $\tilde H$  such that 
$\tilde H \cong H$.
\end{Theorem}

The following is an immediate consequence of Theorem~\ref{VtimesAA} 
for the value $\U=\A^l$:

\begin{Theorem}
\label{SolvableIntoSolvable}
For an arbitrary countable solvable group $H$ of  length $l$
there exists  a two-generated  group $G$ of  length at most $l+2$
(but not, in general, of length $l+1$)
with a subnormal subgroup $\tilde H$  such that 
$\tilde H \cong H$.
\end{Theorem}
\begin{proof}[Proof of Theorem~\ref{VtimesAA}]
We consider the wreath product $H \Wr \Z $. 
For each $h \in H$ we take
a 
$$
   \beta (h) = (\ldots, 1_{H}, 1_{H}, 
   h^{-1}, h^{-1} \ldots)\in H^{(\Zs)}
$$
(the first element $h^{-1}$ is situated 
in the place number $0$). Then
it is easy to compute that
$$
  [1_{\Zs},\beta (h)] = (\ldots,  1_{H}, h, 1_{H},\ldots )
$$
($h$ is in the place number $0$). 
We should now just note that the
 set of all strings on the right 
for all $h \in H$ forms the first copy  of $H$ in $H \Wr \Z$.  

Let
$$
M_1= \langle \Z \cup  \{ \beta (h) | \quad \!\!  h\in H\} \rangle .
$$
$M_1$ is a countable group and 
contains $H$ as a subnormal subgroup. 
Besides an isomorphic copy of $H$ lies in $M_1'$. Therefore we can 
use the point {\bf f} from 
Section~\ref{Themainconstruction} in order 
to build a special element $\overline{m}_1$ in the wreath 
product $M_1 \Wr \Z$ such that the intersection 
$$
\langle -1_{\Zs}, \overline{m}_1 \rangle \cap M_1^{(\Zs)} 
=\overline{M_1'}^{(\Zs)}
$$
contains an isomorphic copy of $H$. Thus 
Theorem~\ref{VtimesAA}
is proved for
the group $
G=\langle -1_{\Zs}, \overline{m}_1 \rangle \in \U \A^2.
$\end{proof}

We note that in these theorems, too, the group $H$ 
is embedded not only
into some appropriate two-generated group $G$ but into the
second derived group of the latter.

Our construction strongly depends on 
the fact that the active groups
of our wreath products contain 
a subgroup isomorphic to $\Z$. Thus
this construction cannot be used in embeddings of, say, finite groups
into finite groups. B.~Neumann and Hanna Neumann prove in 
~\cite{NeumannHNeumann} that 
arbitrary {\sl finite} group 
$H$ of the variety $\U$ is embeddable into
(the second derived group) of 
some {\sl finite} two-generated group 
$H\in \U \A^2$. But the result of Heineken and Lennox 
~\cite{HeinekenLennox} we mentioned
shows that the ``verbal'' 
embeddings of finite groups into finite groups
are in general {\sl not} subnormal.

Finally we observe that the results we proved for
solvable group have no analogs for the cases of abelian 
or nilpotent groups (see the arguments in the introduction).


  \section{On``verbal''embeddings \\ 
  of finite groups into finite groups}
  \label{On``verbal''embeddingsoffinitegroups}

In the sequence we consider embeddings of finite groups only.  Since the
mentioned results of P.~Hall and R.~Dark~\cite{Dark} and of Heineken and 
Lennox~\cite{HeinekenLennox} concern finite {\sl  solvable} groups,
the following Proposition~\ref{FiniteNilpotentIntoFiniteNilpotent} 
shows that those results describe 
in some sense the ``smallest'' possible finite groups  without
finite subnormal ``verbal''  embeddings. 
\begin{Proposition}
\label{FiniteNilpotentIntoFiniteNilpotent}
Let $H$ be  an arbitrary finite nilpotent group and $V$ be an arbitrary
non-trivial word set. Then there exists a finite nilpotent group $G$
with a subnormal subgroup $\tilde H$ such that 
$$
 \tilde H \cong H \quad \mathrm{and} \quad \tilde H \sub V(G).
$$
Moreover if $H$ is a $\Pi$-group for the set $\Pi$ of primes, the 
group $G$ can be chosen a $\Pi$-group, too.
\end{Proposition}
 \begin{proof}
It is clear that if for each group of the family  $\{ H_i \}_{i\in I}$ 
there is  an  embedding  of type  stated in 
Proposition~\ref{FiniteNilpotentIntoFiniteNilpotent} 
(into a group $G_i$, $i\in I$), 
then for
the (direct or cartesian) product $\prod_{i\in I}H_i$ there is an 
embedding of mentioned type as well, namely into the 
(direct or cartesian) product $\prod_{i\in I}G_i$.

And since each finite nilpotent group is a direct product of
its $p$-subgroup, it is sufficient to prove 
Proposition~\ref{FiniteNilpotentIntoFiniteNilpotent} 
for finite $p$-groups.

Modifying   Lemma 1 from~\cite{Heineken} for the case of an 
arbitrary $V$ we obtain for each power $m=p^k$ such a finite $p$-group  
$P$ that its verbal subgroup $V(P)$ is of exponent $m$ and is from the 
center of $P$. Then, as it is shown in the proof of Theorem 3 
in~\cite{Heineken}, the diagonal  
$$
D=\Diag{H \Wr P}= \{ \prod_{t \in P} t^{-1}ht |\quad h \in H  \}
$$
of wreath product $W=H \Wr P$ is contained 
in the verbal subgroup $V(W)$ of $W$. 
If $H$ is a $p$-group $W$ is a $p$-group, too. 
And there is an 
obvious embedding of $H$ into $W$ (onto $D$):
$$
       \nu : h \mapsto \prod_{t \in P} t^{-1}ht.
$$
The base subgroup $H^{(P)}$ of $W$ is normal in $W$ and as a
{\sl nilpotent} group $H^{(P)}$ itself has $D\cong H$ as a
subnormal subgroup.\end{proof}

Are there ``many'' non-nilpotent finite groups $H$ 
for which there is a finite  $G$ such that 
$H \cong \tilde H \sub V(G)$? Are there  ``many'' 
non-solvable finite groups $H$ 
for which one cannot find  a finite  $G$ with this property? The 
following statement shows something more.
\begin{Proposition}
\label{VeryManyGroups}
\begin{enumerate}  
\item For each variety $\W\not=\O$ and an arbitrary
non-trivial word set $V$ there exists an infinite set of 
finite groups $H$ such that $H \notin \W$ and there are finite groups
$G$ and $\tilde{H} \sn G$ such that $H  \cong \tilde{H}   \sub V(G)$;
\item For each variety $\W\not=\O$ there exists an infinite set of 
finite groups $H$ such that $H \notin \W$ and there is no
finite group $G$ with a subgroup $\tilde H  \cong H$  such that 
$\tilde H \sn G'$.
\end{enumerate}
\end{Proposition}
\begin{proof}
1. We have  $\V = \var {\Fin / V(\Fin) } \not= \O$. 
Thus $\W\cdot\V\not=\O$ (see in~\cite{HannaNeumann}).
Since the set of all finite groups generates the variety 
$\O$ of all groups, there is a finite group $G \notin \W\cdot\V$.
Thus $V(G)\notin \W$ because $G$ is an extension of $V(G)$ by 
$G/V(G)$. The finite group $H=V(G)$ satisfies the requirement of
the first point. We can now repeat this step for the variety
$$
\W_1 = \var {\W, H} \not= \O
$$ 
and obtain a new finite group $H_1 \notin \W_1$. The next finite
groups $H_2,H_3,\ldots$ are found by induction.

2. Let $N$ be a finite complete
group that cannot be subnormally 
embedded into the derived group of some finite group. We can 
find such a group according to the result of Heineken and 
Lennox~\cite{HeinekenLennox}. And let $R$ be a finite group with 
and $R\notin \W$.
We can  find  such a group, too, because $\O \not= \W$ and
the set of finite groups generates $\O$.
The group $H= N \times R$ does not belong to $\W$ and does
not have the mentioned subnormal embedding for $N$ does
not have the latter. Taking other values for $R$ we will
obtain  infinitely many examples of group $H$. \end{proof}

Let us observe that  since the set of all varieties 
(i. e. the set of all word sets) is 
of continnum cardinality~\cite{HannaNeumann} and since 
the set of all finite groups
is countable, there exist finite groups that satisfy
the statement 1 or 2 of previous proposition for
$\aleph$ word sets simultaneously. 

We can repeat the argument of Proposition~\ref{VeryManyGroups}  
to obtain finite groups with
needed embedding ``between'' two varieties that can be very 
``near'' one to other.

\begin{Proposition}
\label{ManyGroupsInSmallClasses}
Let $\U$ be an arbitrary non-trivial locally finite
variety and $V$ be a non-trivial word set such that 
$\V = \var {\Fin / V(\Fin) }$ is a locally finite
variety, too. Then for each subvariety $\W$ of $\U$ there
is a finite group in the set $\U \backslash \W$ such that there
exists a finite $G$ with a subnormal subgroup $\tilde H \cong H$, 
$\tilde H \sub V(G)$.
\end{Proposition}
\begin{proof}
Since $\U, \W  \not=\O$ and since multiplication of varieties is a
monotonic operation we have
$$
\W \cdot \V < \U\cdot \V
$$
(see  21.21 in~\cite{HannaNeumann}). The product of two locally finite
varieties is locally finite (O.~Yu.~Schmidt ~\cite{Robinson}). Thus the 
sets ${\mathcal F}_{\W \cdot \V}$ and ${\mathcal F}_{\U \cdot \V}$ of finite groups of 
products $\W \cdot \V $ and $\U \cdot \V$ 
are not equal too. Let 
$1\not= Q \in {\mathcal F}_\U \backslash {\mathcal F}_\W$. Since $Q/V(Q)\in \V$,
necessarily $H=V(Q)\notin \W$. \end{proof} 

We can get similar facts for even finer layers between classes of 
groups~\cite{MikaelianOn"Dividing"Groups}, but it is already 
clear that the finite groups we are looking for are ``many'' and
that the classes of such group are ``unrestricted''. Let us set
the question we are interested in as a {\sl problem}:
\begin{quote}
{\it Let $H$ be a finite group and $V\sub \Fin$ be a non-trivial word set.
Find a criterion under which there exists   a finite 
group $G$ with a subnormal subgroup  $\tilde H$ such that the latter
is isomorphic to $H$ and lies in the verbal subgroup $V(G)$.}
\end{quote}

\vskip6mm

Mathematisches Institut 

Universit\"at W\"urzburg

Am Hubland

D-97074 W\"urzburg 

Germany 

\vskip2mm

e-mail: mikaelian@e-math.ams.org

vahagn.mikaelian@mail.uni-wuerzburg.de

http://www.crosswinds.net/ $\tilde{ }$ mikaelian/home.html

\end{document}